\documentstyle{amsart}
\input{bull-art}
\newcommand{\sect}[1]{\section{#1} \setcounter{equation}{0}}

\begin{document}
\def\currentvolume{27}
\def\currentissue{1}
\def\currentyear{1992}
\def\currentmonth{July}
\def\copyrightyear{1992}
\def\currentpages{159-164}
\title{A counterexample to the Arakelyan Conjecture}
\author{A. Eremenko}
\date{September 11, 1991 and, in revised form,
January 28, 1992}
\address{Alexandre Eremenko, Department of Mathematics, 
Purdue University,
West\break
Lafayette, Indiana 47907 USA}
\email{eremenko@@math.purdue.edu}
\thanks{Supported by NSF grant DMS-9101798}
\maketitle

\begin{abstract}
A ``self--similar'' example is constructed that shows that a
conjecture of N. U. Arakelyan on the order of decrease of 
deficiencies of
an entire function of finite order is not true.
\end{abstract}
\maketitle

\sect{Introduction}

Let $f$ be an entire function and $\delta(a,f)$ denotes the
Nevanlinna deficiency of $f$ at the point $a\in 
\overline{C}$.
Standard references are [8,\ 9]. (No knowedge of Nevanlinna
theory is necessary to understand this paper. We really deal
with a problem of potential theory.) 
Since $\delta(a,f)\geq 0$ and the deficiency relation of 
Nevanlinna 
states that
$$\sum_{a\in \overline {C}}\delta(a,f)\leq 2,$$
it follows that the set of deficient values, that is, 
$\{a:\delta(a,f)>0\}$,
is at most countable. We denote the sequence of 
deficiencies by $\{\delta_n\}$.
In 1966 Arakelyan \cite{A} (see also \cite{GO} or \cite{F})
constructed the first example of
an entire function of {\em finite
order} having infinitely many deficient values. In this 
example the deficiencies
satisfy
\begin{equation}
\label{conjecture}
\sum_{n=1}^{\infty}\frac{1}{\log (1/\delta_n)}<\infty,
\end{equation}
and he conjectured that (\ref{conjecture})
is true for every entire function of finite order.
Another method of constructing such examples was proposed 
in \cite{E1},
but the function in \cite{E1} 
also satisfies
(\ref{conjecture}).

For meromorphic functions of finite order  Weitsman 
\cite{W} proved
\begin{equation}
\label{weitsman}
\sum_{n=1}^{\infty}{\delta_n}^{1/3}<\infty,
\end{equation}
and this is known to be best possible [9,\ 4].
The only known improvement of (\ref{weitsman}) for entire 
functions is
due to Lewis and  Wu \cite{LW}:
\begin{equation}
\label{LW}
\sum_{n=1}^{\infty}\delta_n^{1/3-\epsilon_0}<\infty,
\end{equation}
where $\epsilon_0$ is an absolute constant. In fact, the 
value
$\epsilon_0=2^{-264}$ is given in \cite{LW}.

In this note we will give a construction that produces an 
entire
function of finite order having infinitely many 
deficiencies $\delta_n$
with the property
\begin{equation}
\label{result}
\delta_n\geq c^{-n},
\end{equation}
where $c>1$ is a constant. Thus Arakelyan's conjecture 
(\ref{conjecture})
fails.

Of course a substantial gap still remains between the 
theorem of
Lewis and Wu and
our example. It is natural to ask whether
\begin{equation}
\sum_{n=1}^{\infty}\frac{1}{\log^{1+
\epsilon}(1/\delta_n)}<\infty
\end{equation}
is true with arbitrary $\epsilon>0$ for entire functions 
of finite order.

It is more or less well known that the problem of 
estimating deficiencies
for entire functions of finite order is equivalent to a 
problem of
potential theory. Namely, the following statements are 
equivalent:

{\bf A.} Given any $\rho>1/2$ and a sequence of complex 
numbers $a_n$, there
exists an entire function $f$ of order $\rho$ with the 
property 
$\delta(a_n,f)\geq c\delta_n$ with some constant $c>0$.

{\bf B.} There exist a bounded subharmonic function $u$ in 
the annulus
$A=\{z:1<|z|<2\}$ and
disjoint open sets $E_n\subset A,\;1\leq n<\infty$ with the
following properties:

(i) Each $E_n$ is a union of some components of the set 
$\{z\in A:u(z)<0\}$;

(ii) for every $r\in [1,2]$
$$\int_{\{\theta:re^{i\theta}\in 
E_n\}}u(re^{i\theta})\,d\theta\leq -\delta_n.$$

We indicate briefly how to prove the equivalence.
To prove {\bf A}$\rightarrow${\bf B} we take a sequence of 
P\'{o}lya peaks
\cite[p. 101]{H}
$r_k$ for $\log M(r,f)$ and consider the sequence of 
subharmonic
functions
$$u_k(z)=\frac{\log |f'(r_kz)|}{\log M(r_k,f)},\qquad 
|z|<2.$$
 This sequence
is precompact in an appropriate topology and we may take a 
subsequence that
converges to a subharmonic function $u$. If $f$ has 
deficient values then
$u$ satisfies (i) and (ii). See [1,\ 6] for details.

To prove {\bf B}$\rightarrow${\bf A} we apply the 
construction from \cite{E1}
that involves an extension of $u$ to a subharmonic 
function in
${\bold C}$ with the property of self-similarity: 
$u(2z)=ku(z),\;k=
\mathrm{const}
>0$,
approximation
of $u$ by the logarithm of modulus of an entire function 
$g$ and
performing a quasi-conformal modification on the function 
$g$ that
produces the entire function $f$ satisfying {\bf A}. It is 
also plausible
that Arakelyan's original method could be applied 
directly as soon as a subharmonic
function with the properties (i) and (ii) is constructed.

{\em Remark}. The above-mentioned paper of Lewis and Wu 
contains also
the solution of a problem of Littlewood on the upper 
estimate of 
mean spherical derivative of a polynomial. The connection 
between the two
problems seems somewhat obscure. An example that gives a 
lower estimate in
the Littlewood's problem was constructed in \cite{E3} 
using some self-similar
sets arising in the iteration theory of polynomials. It is 
interesting 
that the example we are going to construct now also has 
the property of
self-similarity. Instead of iteration of a polynomial here 
the crucial
role is played by a semigroup of M\"{o}bius 
transformations of the plane. 

\section{The example}

Consider the semigroup $\Gamma$ generated by $z\mapsto 
z\pm 1$ and
$z\mapsto z/2$. We have
$$\Gamma=\{\gamma_{n,k}:n=0,1,2,\ldots ;\,k=0,\pm 1,\pm 
2,\ldots\},$$
where $\gamma_{n,k}(z)=2^{-n}(z+k).$

Denote by $S^+_{0,0}$ the square
$$S^+_{0,0}=\{z:|\Re z|\leq\tfrac{3}{10},\;|\Im 
z-1|\leq\tfrac{3}{10}\}$$
and set $S^+_{n,k}={\gamma}_{n,k}(S^+_{0,0})$. It is easy 
to see that
the squares $S^+_{n,k}$ are disjoint. Consider
the domain
$$D_0=\{z:0<\Im z<\tfrac{4}{3}\}\backslash \bigcup_{n,k}S^+
_{n,k}.$$
The boundary $\partial D_0$ consists of the real axis, 
boundaries of
the squares and horizontal line $l_0=\{z:\Im z=4/3\}$. The 
domain
$D_0$ is $\Gamma$-invariant 
and the transformation $z\mapsto z/2$ 
maps $D_0$ onto
$$D_1=\{z:0<\Im z<\tfrac{2}{3}\}\backslash\bigcup_{n,k}S^+
_{n,k}\subset D_0.$$
The boundary of $D_1$ consists of the real axis, 
boundaries of some squares,
and the horizontal line $l_1=\{z:\Im z=2/3\}\subset D_0$.

Let $u$ be the harmonic function in $D_0$ that solves the 
Dirichlet problem
$$u(z)=1,\qquad z\in l_0,$$
$$u(z)=0,\qquad z\in\partial D_0\backslash l_0.$$
This Dirichlet problem has a unique solution. So we 
conclude from translation
invariance that
\begin{equation}
\label{periodicity}
u(z+1)=u(z),\qquad z\in D_0.
\end{equation}
It follows that the function $u$ has a positive minimum 
$M^{-1}<1$ on
the line $l_1\subset D_0$. Comparing $u(z)$ and $u(2z)$ on 
$\partial D_1$
and using the maximum principle, we conclude that 
$u(2z)\leq Mu(z),\;z\in D_1$,
which is equivalent to
\begin{equation}
\label{intermediate}
u(z)\leq Mu(z/2),\qquad z\in D_0.
\end{equation}
It follows from (\ref{periodicity}) and 
(\ref{intermediate}) that
\begin{equation}
\label{selfsimilarity}
u(\gamma_{n,k}(z))\geq M^{-n}u(z),\qquad z\in D_0.
\end{equation}

Now we are going to extend $u$ to the strip
$$S^+=\{z:0<\Im z <\tfrac{4}{3}\},$$
that is, to define $u$ in the squares.
We start by defining $u$ in $S^+_{0,0}$. The normal 
derivative
(in the direction of the outward normal to the boundary of 
the square)
of $u$ has positive infimum on $\partial S^+_{0,0}$; it 
tends to $+\infty$ as we approach a corner of the square. 
Denote by $G>0$
the Green 
function for $S^+_{0,0}$ with the pole at the point 
$i=\sqrt{-1}$. It is
clear that the normal derivative of $G$ on the boundary of 
the square
is bounded (it tends to zero as we approach a corner). Set
$$u(z)=-tG(z),\qquad z\in S^+_{0,0},$$
where $t>0$. If $t$ is small enough we obtain a 
subharmonic extension
of $u$ into $S^+_{0,0}$, because the jump of the normal 
derivative
will be positive as we cross the boundary of the square 
from inside.
Fix such $t$, and extend $u$ to the remaining squares by 
the formula
\begin{equation}
\label{extension}
u(\gamma_{n,k}(z))=M^{-n}u(z),\qquad z\in S^+_{0,0}.
\end{equation}
It follows from (\ref{selfsimilarity})
that the normal derivative always has a positive jump as 
we cross the boundary
of $S^+_{n,k}$, so the extended function is subharmonic in 
$S^+$.

Now consider the smaller squares
$$K^+_{0,0}=\{z:|\Re z|\leq\tfrac{2}{7},\;|\Im 
z-1|\leq\tfrac{2}{7}\}\subset
S^+_{0,0},\qquad
K^+_{n,k}=\gamma_{n,k}(K^+_{0,0})\subset S^+_{n,k}.$$
It follows from (\ref{extension}) that
\begin{equation}
\label{negative}
u(z)\leq-\beta M^{-n},\qquad z\in K^+_{n,k}
\end{equation}
for some $\beta >0$ and all $n$ and $k$.

Now we are going to extend $u$ to the strip 
$S=\overline{S^+\cup S^-}$ where
$$S^-=\{z:-\tfrac{4}{3}<\Im z <0\}.$$
To do this we repeat the above construction starting with 
the square
$$S^-_{0,0}=\{z:|\Re 
z-\tfrac{1}{2}|\leq\tfrac{3}{10},\;|\Im z+1|
\leq\tfrac{3}{10}\}$$
and using the same semigroup $\Gamma$. We obtain the squares
$$S^-_{n,k}=\gamma_{n,k}(S^-_{0,0})\quad \text{and}\quad
K^-_{n,k}=\gamma_{n,k}(K^-_{0,0}),$$
where
$$K^-_{0,0}=\{z:|\Re z-\tfrac{1}{2}|\leq 
\tfrac{2}{7},\;|\Im z+1|\leq
\tfrac{2}{7}\},$$
and the function $u_1$ subharmonic in $S^-$ that satisfies 
the inequality
similar to (\ref{negative}):
\begin{equation}
\label{u1negative}
u_1(z)\leq-\beta_1 M_1^{-n},\qquad z\in K^-_{n,k}
\end{equation}
with some $\beta_1>0$ and $M_1>1$.

Extend $u$ to $S=\overline{S^+\cup S^-}$ by setting 
$u(z)=u_1(z),\;z\in S^-$
and $u(x)=0,\;x\in {\bold R}$. The extended
function $u$ is continuous in $S$. We will prove that it 
is subharmonic in $S$.

Consider the strips
$\Pi_n=\{z:|\Im z|<\tfrac{4}{3} 2^{-n}\}.$
Define the functions $v_n$ in the following way:
$v_n(z)=u(z),\;z\in S\backslash \Pi_n;\;v_n$ are continuous 
in $S$ and harmonic in $\Pi_n$. Then $v_n(x)>0,\;x\in 
{\bold R}$, and
it follows from the maximum principle
(applied to $\Pi^+_n=\{z:0<\Im z<\frac{4}{3} 2^{-n}\}$)
that $v_n\geq u$ in $S$. We conclude that $v_n$ are 
subharmonic because
the sub-mean value property holds in every point of $S$. 
Furthermore,
it is evident that $v_n\rightarrow u$ uniformly in $S$ as 
$n\rightarrow
\infty$, so $u$ is subharmonic in $S$.

Denote $K^+_n=\bigcup_k K^+_{n,k}$ and $K^-_n=\bigcup_k 
K^-_{n,k}$
and remark that each vertical line $\Re z=x_0$ intersects 
$K^+_n\cup K^-_n$.
Indeed the projection of $K^+_n$ onto the real axis is
$$\bigcup_{k\in {\bold Z}}\{x:|x-2^{-n}k|\leq\tfrac{2}{7} 
2^{-n}\}$$
and the projection of $K^-_n$ is
$$\bigcup_{k\in {\bold Z}}\{x:|x-2^{-n}(k-\tfrac{1}{2})|\leq
\tfrac{2}{7} 2^{-n}\}.$$
It is clear that the union of these two sets is the whole 
real axis.

Now set $E_n=\bigcup_k(S^+_{n,k}\cup S^-_{n,k})$ and $K_n=
K^+_n\cup K^-_n$. Then each $E_n$ is a union of some 
components of the set
$\{z\in S:u(z)<0\}$ and for every vertical line
$l=\{z:\Re z=x_0\}$ the length of intersection 
$l\cap K_n$ is at least $b^{-n}$ for some $b>1$.
So we have  
in view of (\ref{negative}) and
(\ref{u1negative}):
$$\int_{l\cap E_n} u(x_0+iy)\,dy\leq
\int_{l\cap K_n} u(x_0+iy)\,dy\leq -c^{-n}$$
with some constant $c>1$.

It remains to make a change of variable 
$z=\frac{4}{3\varepsilon}\log \zeta,
\;\zeta\in Q=\{\zeta:1<|\zeta|<2,\;|\arg 
\zeta|<\varepsilon\}$, and to
extend the function $u(z(\zeta))$ (in any desired manner)
from $Q$ to a subharmonic function in
the annulus $\{1<|\zeta|<2\}$. The extended function will 
have the
properties
 (i) and (ii) of {\bf A} with $\delta_n\geq c^{-n},\;c>1.$

\end{document}